\documentstyle[11pt]{article}
\begin{document}

\vspace*{1in}
\centerline {{\Large {\bf Dynamics of Certain Smooth
One-dimensional Mappings}}}

\vskip20pt
\centerline {\Large \bf III. Scaling function geometry}

\vskip30pt
\large
\centerline{Yunping Jiang }
\centerline{Institute for Mathematical Sciences, SUNY at Stony Brook}
\centerline{Stony Brook, L.I., NY 11794}
\vskip20pt
\centerline{December 1, 1990}

\vskip40pt
\centerline{\em Dedicated to Professor Shantao Liao on the occasion of his seventieth birthday}

\vskip50pt

\centerline{ {\Large \bf Abstract}}

\vskip20pt
We study scaling function geometry. We show the existence of the scaling
function of a geometrically finite one-dimensional mapping. This scaling
function is discontinuous. We prove that the
scaling
function and the asymmetries at the critical points of a
geometrically finite one-dimensional mapping
form a complete set of
$C^{1}$-invariants within
a topological conjugacy class.

\pagebreak
\centerline{{\Large \bf Contents}}

\vskip20pt
\noindent \S 1 Introduction.

\vskip10pt
\noindent \S 2 The Scaling Structure of a Markov Mapping.

\vskip5pt
\S 2.1 A Markov partition.

\vskip5pt
\S 2.2 The symbolic and the dual symbolic spaces.

\vskip5pt
\S 2.3	 The signed scaling function.

\vskip5pt
\S 2.4	Some properties of a signed scaling function.

\vskip10pt
\noindent \S 3 The
Scaling Function of a Geometrically Finite
One-dimensional

Mapping.

\vskip5pt
\S 3.1 The existence of the signed scaling function.

\vskip5pt
\S 3.2 The $C^{1+\alpha}$-classification.

\pagebreak

\centerline{ \Large \bf \S 1 Introduction}

\vskip10pt

{\bf Smooth classification.} Two smooth mappings $f$ and $g$ from a
one-dimensional manifold $M$ to itself are
topologically conjugate if there is a homeomorphism $h$ from $M$ to
itself such that $f\circ h =h\circ g$. The homeomorphism $h$
is not usually a smooth diffeomorphism;
if is is, then all the eigenvalues of $f$ and $g$ at the
corresponding periodic points have to be the same.
We say $f$ and $g$
are \underline{smoothly conjugate}
if the homeomorphism $h$ is a smooth diffeomorphism,

\vskip5pt
A celebrated theorem first proved by M. Herman [H] says that
any circle diffeomorphism with a diophantine
rotation
number is smoothly conjugate to the rigid rotation by that number. In
this case, the circle diffeomorphism lacks periodic points. Thus the
topological invariant, the rotation number, is the complete smooth
invariant.

\vskip5pt
A theorem proved by M. Shub and D. Sullivan [SS] says that any two
smooth
orientation-preserving expanding endomorphisms of the circle are
smoothly
conjugate if they are topologically conjugate and the conjugacy
is absolutely continuous.
D. Sullivan
[S1] also showed that the set of eigenvalues at periodic
points of a
smooth orientation-preserving expanding endomorphism of
the circle forms
a complete set of smooth invariants within a topological conjugacy
class.

\vskip5pt
A recent theorem proved by D. L. Llave and R. Moriy\'on [LM] says that
any two
Anosov diffeomorphisms on the torus are smoothly conjugate if they are
topologically conjugate and
have the same eigenvalues at
corresponding
periodic points.  Thus the set of eigenvalues at periodic points of
an Anosov diffeomorphism on the torus forms a complete set of smooth
invariants within a topological conjugacy class.

\vskip5pt
A recent work in [J1], I provided a smooth classification of the space
of generalized Ulam-von Neumann transformations.  These are
certain smooth interval mappings topologically conjugate to the
mapping $q(x)=-x^{2} +2$ of the interval $[-2, 2]$.
I classified
this space up to smooth
equivalence by showing that
all the eigenvalues at periodic points,
the type of power law
at the critical point, and a quantity which we call the asymmetry at
the critical point form a
complete and optimal set of smooth invariants.

\vskip5pt
{\bf What we would like to say in this paper.}
We study continuously geometrically finite
one-dimensional mappings (see [J3]). These are a subspace of
$C^{1+\alpha }$ one-dimensional mappings with finitely many,
critically finite power law critical points. We
concentrate on scaling function geometry. We show the
existence of the scaling
function of a geometrically finite one-dimensional mapping. We
study the rigidity on the space of geometrically
finite one-dimensional mappings.

\vskip5pt
Suppose $M$ is an oriented connected compact one-dimensional
$C^{2}$-Riemannian manifold with Riemannian metric $dx^{2}$ and
associated length element $dx$. Suppose $f$ is
a mapping from $M$ into $M$.
A \underline{Markov partition} $\eta_{1}$ of $M$ by $f$ is a set $\{
I_{0}$, $\cdots $,
$I_{m}\}$ of closed intervals of $M$ such that $(a)$ $I_{0}$, $\cdots $,
$I_{m}$ have pairwise disjoint
interiors, $(b)$ the union $\cup_{i=0}^{m}I_{i}$ of the intervals is
$M$,
$(c)$ the restriction of $f$ to each interval $I_{i}$ is injective and
continuous, and $(d)$ the image of $I_{i}$ under $f$ is the union of some
intervals in $\eta_{1}$. We say $f$ is a \underline{Markov mapping} if
it has a Markov partition $\eta_{1}$ of $M$ by $f$.

\vskip5pt
\newcommand{\sn}{w_{n}=r_{i_{0}}\cdots r_{i_{n}}}
\newcommand{\gs}{g_{r_{i_{0}}\cdots r_{i_{n}}}}
\newcommand{\gls}{g_{r_{i_{0}}}\circ \cdots \circ g_{r_{i_{n}}}}
\newcommand{\Ig}{I_{r_{i_{0}} \cdots r_{i_{n}}}}

Suppose $f$ is a Markov mapping and $\eta_{1}$ is a fixed Markov
partition of $M$ by $f$.
We
use the symbols $0$, $\cdots$, $m$ to name the intervals in the
partition $\eta_{1}$ and use $g_{i}$ to denote the inverse of the
restriction of $f$
to the interval with name $i$.	Let symbol $r_{i}$ be $+i$ if
$g_{i}$ is
orientation-preserving and be $-i$ if $g_{i}$ is orientation-reversing.
Suppose $\sn$ is a sequence of the symbols $\{ r_{0}$, $\cdots$,
$r_{m}\}$.
We say it is a \underline{suitable sequence} of length $n+1$ if
$I_{i_{k}}\subset
f(I_{i_{k-1}})$ for any $k=0$, $\cdots$, $n$.  In the other words, it is
suitable if the interval $I_{i_{k}}$ is in the domain of $g_{i_{k-1}}$ for
any $k=0$, $\cdots$, $n$.
For a suitable sequence $w_{n}=r_{i_{0}}\cdots r_{i_{n}}$,
we use
$g_{w_{n}}$ to denote the composition $g_{i_{0}}\circ \cdots \circ
g_{i_{n}}$ and use $I_{w_{n}}$ to denote
the image of $f(I_{i_{n}})$ under $g_{w_{n}}$. We call $w_{n}$ the
name of the interval $I_{w_{n}}$.
We may read the name $w_{n}$ either from the left to the right or from
the right to the left.

\vskip5pt
Suppose we read all the names from the left to the right. Then we
get the set $\Sigma_{f}=\{ a= r_{i_{0}}r_{i_{1}}\cdots \}$ of
infinite suitable sequences which start from the left and extend to the
right. Suppose $\sigma_{f}: \Sigma_{f} \mapsto
\Sigma_{f}$ is the
shift mapping which knocks off the first symbol in the left of an
infinite suitable sequence in $\Sigma_{f}$.
The space $(\Sigma_{f} , \sigma_{f})$ is the phase space of the Markov
mapping $f$.

\vskip5pt
Let us now read all the names from the right to the left.	We
then get the set $\Sigma^{*}_{f}=\{ a^{*}= \cdots r_{i_{1}}r_{i_{0}}\}$
of
infinite suitable sequences which start from the right and extend to the
left. Suppose $\sigma_{f}^{*}: \Sigma_{f}^{*} \mapsto
\Sigma_{f}^{*}$ is the
shift mapping which knocks off the first symbol in the right of an
infinite suitable sequence in $\Sigma_{f}^{*}$.
We call the space $(\Sigma_{f}^{*}, \sigma_{f}^{*} )$ the
\underline{dual} \underline{space} of the Markov mapping $f$.

\vskip5pt
The mapping $sign: \eta_{n}\mapsto \{ -1,$
$1\}$ is defined
by $sign (I_{w_{n}})$ where $sign (I_{w_{n}})$ is $1$ if the number of
$-$ in the sequence $w_{n}$ is even and $sign(I_{w_{n}})$ is $-1$ if
the number of $-$ in the sequence $w_{n}$ is odd.
For an infinite suitable sequence $a^{*}= \cdots w_{n}$ in
$\Sigma_{f}^{*}$, let $\sigma_{f}^{*}(a^{*}) = \cdots v_{n-1}$ where
$w_{n}=v_{n-1}r_{i_{0}}$.  We use
$s(w_{n})$ to denote the ratio
\[ \frac{sign(I_{w_{n}})|I_{w_{n}}|}{sign(I_{v_{n-1}})|I_{v_{n-1}}|}\]
and call it the signed scale at $w_{n}=v_{n-1}r_{i_{0}}$. We also call the
absolute value of the signed scale the scale. If the limit
$s_{f}(a^{*})$ of the sequence $\{
s(w_{n})\}_{n=0}^{\infty}$ of the signed scales exists as $n$ goes to
infinity, then we say there is the signed scale
at $a^{*}$.
If there is the signed scale at every point in $\Sigma_{f}^{*}$, then we
define a
function $s_{f}: \Sigma_{f}^{*}\mapsto {\bf R^{1}}$ as $s_{f}(a^{*})$.
We call this function the \underline{signed scaling function} of $f$ and
its absolute value the \underline{scaling function} of $f$.

\vskip5pt
The scaling function was first defined by M. Feigenbaum [F] to describe the
universal geometric structure of the attractors obtained by period
doubling.
D. Sullivan [S2] defined in mathematics the scaling function for a Cantor set
which is the maximal invariant set of a $C^{1+\alpha}$-expanding mapping for
some $0<\alpha \leq 1$. He gave a complete $C^{1+\alpha}$-classification of
these Cantor sets on the line by their scaling functions and used this
classification in the study of the universal geometric structure of the
attractors obtain by period doubling.
The definition of a signed scaling function in
this
paper (also see [J4]) generalizes Sullivan's idea to a Markov mapping.

\vskip5pt
We show some basic properties of the signed scaling function of a Markov
mapping in \S 2.4 (Proposition 1 to Proposition 4).

\vskip5pt
A geometrically finite one-dimensional mapping $f$ has been defined in
[J3], which is a certain Markov mapping with finitely many, critically finite
power law critical points (see \S 3 for a definition).
The fixed Markov partition $\eta_{1}$ of a geometrically finite
one-dimensional mapping $f$ is the set of the closures of the
intervals of the complement of the critical orbits of $f$.
One of the main theorems in this paper is the
following:

\vskip5pt
{\sc Theorem A.} {\em
Suppose $f$ is a geometrically finite one-dimensional mapping.
Then there is the signed scaling function $s_{f}:
\Sigma^{*}_{f}\mapsto {\bf R}^{1}$ of $f$.}

\vskip5pt
The proof of this theorem is an application of the
$C^{1+\alpha}$-Denjoy-Koebe distortion lemma in [J2] (see also [J3]).

\vskip5pt
Suppose $f$ is a geometrically finite one-dimensional mapping.
A critical point $c$ of $f$ is a power law critical point of $f$
if there is a number
$\gamma > 1$ such that the limits of $f'(x)/|x-c|^{\gamma -1}$ exist
and equal nonzero numbers as $x$ goes to $c$ from below and from above.
We call the number $\gamma$ the exponent of $f$ at the power law
critical point $c$.
Two corollaries of Theorem A are the following:

\vskip5pt
{\sc Corollary A1.} {\em
Suppose $f$ is a geometrically finite one-dimensional mapping.
Then the scaling function $|s_{f}|:
\Sigma_{f}^{*}\mapsto {\bf R}^{1}$ of $f$ is
discontinuous.}

\vskip5pt
{\sc Corollary A2.} {\em
Suppose $f$ is a geometrically finite one-dimensional mapping.
Then the exponent $\gamma$ of $f$ at a power law critical point $c$
can be calculated by the scaling function
$s_{f}:\Sigma_{f}^{*}\mapsto {\bf R}^{1}$.}

\vskip5pt
Suppose $f$ is a geometrically finite one-dimensional mapping.
We say an object is a \underline{$C^{1}$-invariant} of $f$ if it is
the same for $f$ and for $h\circ f\circ h^{-1}$ whenever $h$ is an
orientation-preserving $C^{1}$-diffeomorphism.
The \underline{asymmetry} of $f$ at a power law critical point $c$ of
$f$
is the limit of $f'(x)/f'(-x+2c)$ as $x$ tends to $c$ from below. It is
a $C^{1}$-invariant of $f$ (see [J1]).
The signed scaling function $s_{f}$ of $f$
is a $C^{1}$-invariant too (see Proposition 1 in \S 2.4).
Another main result in this paper is that the scaling function
$s_{f}$ of $f$ and the asymmetries
of $f$ at all the critical points of $f$ form a complete set of
$C^{1}$-invariants
within a topological conjugacy class as follows.

\vskip5pt
{\sc Theorem B.} {\em Suppose $f$ and $g$ are geometrically finite and
topologically conjugate by an orientation-preserving homeomorphism $h$.
Then
$f$ and $g$ are $C^{1}$-conjugate if and only if
they have the same signed scaling function and the same asymmetries at
the corresponding critical points.}

\vskip5pt
Actually, we can say
more on the smoothness of the conjugating mapping $h$ as follows.

\vskip5pt
{\sc Corollary B1.}
{\em
Suppose $f$ and $g$ are $C^{1+\alpha}$-geometrically finite
one-dimensional mappings for some $0< \alpha \leq 1$.
Furthermore suppose they are topologically conjugate by an
orientation-preserving homeomorphism $h$.
If $f$ and $g$ have the same signed scaling function and
the same asymmetries at the corresponding critical points, then
$h$ is a $C^{1+\alpha
}$-diffeomorphism. }

\vskip5pt
Suppose $f:M\mapsto M$ is a geometrically finite one-dimensional mapping
and
$s_{f}:\Sigma_{f}\mapsto \Sigma_{f}$ is the signed scaling function of
$f$. The eigenvalue $e_{f}(p)=(f^{\circ n})'(p)$
of $f$ at a periodic point $p$ of period $n$ and
the exponent $\gamma$ of $f$ at a critical point $c$
can be calculated by the signed scaling function $s_{f}$ of $f$ showed
by Proposition 2 and
by Corollary A2.  Both of them are then clearly
$C^{1}$-invariants. Moreover, in the case that
the set of periodic points of $f$ is dense in $M$,
we show that the
eigenvalues of $f$ at
periodic points and the exponents and the asymmetries of $f$ at critical
points form a complete $C^{1}$-invariants within a topologically
conjugate class as
follows.

\vskip5pt
{\sc Theorem C.} {\em Suppose $f$ and $g$ are
$C^{1+\alpha}$-geometrically finite one-dimensional mappings for some
$0< \alpha \leq 1$. Furthermore, suppose $f$ and $g$ are topologically
conjugate by an orientation-preserving homeomorphism $h$ and suppose the
set of periodic points of $f$ is dense in $M$. If $f$ and $g$ have the
same eigenvalues at the corresponding
periodic points and the same exponents at the
corresponding critical
points, then they have the same scaling functions.
Moreover, if $f$ and $g$ have also the same asymmetries at the
corresponding critical points, then $h$ is a $C^{1+\alpha
}$-diffeomorphism.}

\vskip5pt
In the case that the set of periodic points of $f$ is not dense in $M$,
it seems that 
the eigenvalues of $f$ at
periodic points and the exponents and the asymmetries of $f$ at critical
points are not enough to form a complete $C^{1}$-invariants within a topologically
conjugate class for there may be a Cantor set in $M$ which is an invariant
set of some iterate of $f$ (see the results in [CP]). But we are still
interested in the following question.

\vskip5pt
{\sc Question 1.} {\em Suppose $f: M\mapsto M$ is a geometrically finite
one-dimensional mapping. Do the eigenvalues of $f$
at the periodic points, the exponents and the asymmetries of $f$ at the
critical points of $f$
form a complete set of $C^{1}$-invariants within a topologically conjugate
class ?}

\vskip20pt
{\bf Acknowledgment.} The author would like to
thank Dennis Sullivan and John Milnor for their constant encouragement
and many useful conversations. I would
also like to thank Benjamin Bielefeld, Karen Brucks, Elise Cawley, Mikhail Lyubich,
Scott Sutherland, Grzegorz Swiatek,
Folkert
Tangerman and Peter Veerman
for many useful conversations and help.

\vskip10pt
\centerline {\Large \bf \S 2 Scaling Structure of a Markov Mapping}

\vskip5pt
Suppose $M$ is an oriented connected compact one-dimensional
$C^{2}$-Riemannian manifold with Riemannian metric $dx^{2}$ and
associated length element $dx$. Suppose $f: M\mapsto M$ is a
piecewise continuous mapping.

\vskip5pt
\noindent {\bf \S 2.1 A Markov partition.}

\vskip3pt
A Markov partition $\eta_{1}$ of $M$ by $f$ is a set $\{ I_{0}$, $\cdots
$,
$I_{m}\}$ of closed intervals of $M$ such that $(a)$ $I_{0}$, $\cdots $,
$I_{m}$ have pairwise disjoint
interiors, $(b)$ the union $\cup_{i=0}^{m}I_{i}$ of the intervals is
$M$,
$(c)$ the restriction of $f$ to each interval $I_{i}$ is injective and
continuous, and $(d)$ the image of $I_{i}$ under $f$ is the union of some
intervals in $\eta_{1}$.

\vskip5pt
{\sc Definition 1.} {\em The mapping $f:M\mapsto M$ is a Markov mapping if
there is a Markov partition of $M$ by $f$.}

\vskip5pt
\noindent {\bf \S 2.2 The symbolic and the dual symbolic spaces.}

\vskip3pt
Suppose $f$ is a Markov mapping and $\eta_{1}$ is a fixed Markov partition
of $M$ by $f$. We use $0$, $\cdots$, $m$
to name the
intervals in $\eta_{1}$ and use $g_{i}$ to denote the inverse of the
restriction of $f$ to the interval with name $i$. Let symbol $r_{i}$ be
$+i$ if $g_{i}$
is orientation-preserving and be $-i$
if $g_{i}$ is orientation-reversing.

\vskip5pt
Suppose $\sn$ is a sequence of the symbols $\{ r_{0}$, $\cdots$,
$r_{m}\}$.
We say it is a \underline{suitable sequence} of length $n+1$ if
$I_{i_{k}}\subset
f(I_{i_{k-1}})$ for any $k=0$, $\cdots$, $n$.  In the other words, it is
suitable if the interval $I_{i_{k}}$ is in the domain of $g_{i_{k-1}}$ for
any $k=0$, $\cdots$, $n$.
Suppose $\sn$ is a suitable sequence. Let $g_{w_{n}} = \gls$
be the composition of $g_{i_{0}}$ to $g_{i_{n}}$ and
let $I_{w_{n}} = g_{w_{n}} (f(I_{i_{n}}))$ be the image of
$f(I_{i_{n}})$ under $g_{w_{n}}$. We call $w_{n}$ the
name of the interval $I_{w_{n}}$.  Suppose $\eta_{n}$ is the set of
the intervals $I_{w_{n}}$ for all suitable sequences $w_{n}$ of
length $n$.  This set is
also a Markov partition of $M$ by $f$. We call it the $n^{th}$-partition of
$M$ by $f$. Let $\lambda_{n}$ be the maximum of the lengths of the
intervals in $\eta_{n}$. \underline{We always assume that $\lambda_{n}$ tends to
zero as $n$ goes to infinity.}

\vskip5pt
Suppose $\Gamma_{n}$ is
the set of all the names $w_{n}$ of the intervals in $\eta_{n}$.
We define a
$(n,k)$-left cylinder for $0\leq k\leq n$ as
\[ [w_{n}^{0}]=\{ w_{n}=r_{i_{0}}\cdots r_{i_{n}}
| \in \Gamma_{n}, r_{i_{l}}=r_{i_{l}}^{0}, \hskip5pt for \hskip5pt l=0,
\cdots ,k \}\] where
$w_{n}^{0}=r_{i_{0}}^{0}\cdots r_{i_{n}}^{0}$ is a fixed suitable
sequence in $\Gamma_{n}$.
All the $(n,k)$-left
cylinders form a
topological basis of $\Gamma_{n}$. We still use $\Gamma_{n}$ to
denote the set $\Gamma_{n}$ with this topological basis.  The
sequence $\{ \Gamma_{n}\}_{n=0}^{\infty}$ of the topological spaces
$\Gamma_{n}$ with the inclusions $I_{n}: \Gamma_{n}\mapsto \Gamma_{n-1}$
forms an inverse limit set. Its inverse limit $\Sigma_{f}=\{
a=r_{i_{0}}r_{i_{1}}\cdots \}$ with the shift mapping $\sigma_{f}:
\Sigma_{f}\mapsto \Sigma_{f}$ which is defined as
$\sigma_{f}(r_{i_{0}}r_{i_{1}}\cdots )=r_{i_{1}}\cdots$ is the phase space of the
dynamical system $f:M\mapsto M$ as follows:

\vskip5pt
{\sc Lemma 1.} {\em There is a continuous mapping $h$
from $\Sigma_{f}$ onto $M$ such that
\[ f\circ h=h\circ \sigma_{f}\]
and the fiber $h^{-1}(x)$ contains at most two points for every $x\in
M$.}

\vskip5pt
We now consider a
$(n,k)$-right cylinder for $0\leq k\leq n$ as
\[ [w_{n}^{0}]=\{
w_{n}=r_{i_{n}}\cdots r_{i_{0}} |\in \Gamma_{n},
r_{i_{l}}=r_{i_{l}}^{0} \hskip5pt for \hskip5pt l=0, \cdots ,k \}\]
where
$w_{n}^{0}=r_{i_{n}}^{0}\cdots r_{i_{0}}^{0}$ is a fixed suitable
sequence in $\Gamma_{n}$.
All the $(n, k)$-right
cylinders form another
topological basis of $\Gamma_{n}$. Let $\Gamma^{*}_{n}$ be the set
$\Gamma_{n}$ with this topological basis. The sequence $\{
\Gamma_{n}^{*} \}_{n=0}^{\infty}$ of the topological spaces
$\Gamma^{*}_{n}$ with the inclusions $I^{*}_{n}: \Gamma^{*}_{n}\mapsto
\Gamma^{*}_{n-1}$
forms an inverse limit set. Its inverse limit $\Sigma^{*}_{f}=\{
a^{*}=\cdots r_{i_{1}}r_{i_{0}} \}$ with the shift mapping
$\sigma^{*}_{f}: \Sigma^{*}_{f}\mapsto \Sigma^{*}_{f}$ which is defined
as $\sigma^{*}_{f}(\cdots r_{i_{1}}r_{i_{0}})=\cdots r_{i_{1}}$ is not the
phase space of the dynamical system $f:M\mapsto M$ any more. We call it the dual
space of $f$. The scaling invariants will be defined on this dual space
as follows.

\vskip5pt
Let $sign: \eta_{n}\mapsto \{ -1,$
$1\}$ be the mapping defined
by $sign (I_{w_{n}})$ where $sign (I_{w_{n}})$ is $1$ if the number of
$-$ in the sequence $w_{n}$ is even and $sign(I_{w_{n}})$ is $-1$ if
the number of $-$ in the sequence $w_{n}$ is odd.
Suppose $a^{*}= \cdots w_{n}$ is a point in
$\Sigma_{f}^{*}$ and $\sigma_{f}^{*}(a^{*})=\cdots v_{n-1}$ where
$w_{n}=v_{n-1}r_{i_{0}}$. We define the signed scale at $w_{n}$
as
\[
s(w_{n})=\frac{sign(I_{w_{n}})|I_{w_{n}}|}{sign(I_{v_{n-1}})|I_{v_{n-1}}|}.\]
We call the
absolute value $|s(w_{n})|$ of the signed scale $s(w_{n})$ at $w_{n}$
the scale at
$w_{n}$. If the limit
\[ s_{f}(a^{*})=\lim_{n\mapsto \infty} s(w_{n})\]
exists, then we say
there is the signed scale
at $a^{*}$. We call the absolute value $|s_{f}(a^{*})|$ of the signed
scale $s(a^{*})$ at $a^{*}$ the scale at $a^{*}$.

\vskip5pt
{\sc Definition 2.} {\em Suppose there is the signed scale at every
point in $\Sigma_{f}^{*}$. Then we call the
function $s_{f}: \Sigma_{f}^{*}\mapsto {\bf R^{1}}$ defined as
the signed scale $s_{f}(a^{*})$
the signed scaling function of $f$ and call its
absolute value $|s_{f}|$ the scaling function of $f$.}

\vskip5pt
\noindent {\bf \S 2.4 Some properties of a signed scaling function.}

\vskip3pt
We show some properties of a signed scaling function (if it exists) of a
Markov
mapping

\vskip5pt
{\sc Definition 3.} {\em Suppose $f:M\mapsto M$ is a continuous mapping.
We
say an object associated with $f$ is a $C^{1}$-invariant of $f$ if
it is the same for $f$ and for $h\circ f\circ h^{-1}$ whenever $h$ is an
orientation-preserving $C^{1}$-diffeomorphism.}

\vskip5pt
The following
proposition
is immediately from the definition of a signed scaling function.

\vskip5pt
{\sc Proposition 1.} {\em Suppose $f: M\mapsto M$ is a Markov
mapping. Then the signed scaling function
$s_{f}:\Sigma_{f}^{*}\mapsto {\bf R}^{1}$ (if it exists)
is a $C^{1}$-invariant of
$f$.}

\vskip5pt
Suppose $f$ is a Markov mapping and $(\sigma_{f}^{*}, \Sigma_{f}^{*})$
is the dual space of $f$.
Let
$P(\sigma_{f}^{*})$ and $P(f)$ be the
sets of the periodic points of $\sigma_{f}^{*}$ and $f$,
respectively.

\vskip5pt
{\sc Proposition 2.} {\em There is a surjective mapping $Q: P(\sigma_{f}^{*}) \mapsto
P(f)$ such that every fiber $Q^{-1}(p)$ contains at most two points.}

\vskip5pt
{\it Proof.} Suppose $a^{*}$ is a point in $P(\sigma_{f}^{*})$. It can
be written in $a^{*}= (w_{n})^{\infty}= \cdots w_{k}w_{k}$
where
$w_{k}=r_{i_{k-1}}\cdots r_{i_{0}}$ is a finite sequence. The intervals
with the names $(w_{k})^{l}$ satisfy that
\[ \cdots \subset I_{(w_{k})^{l}} \subset
I_{(w_{k})^{l-1}} \subset \cdots \subset I_{w_{k}} \subset
I_{r_{i_{k-1}}}.\]
Let
$\{ p \} = \cap_{l=0}^{\infty} I_{(w_{n})^{l}}$, it is easy to check that
$ f^{\circ k}(p) = p. $
We define $Q: P(\sigma_{f}^{*}) \mapsto P(f)$
as
$ Q((w_{k})^{\infty}) = p $ where $p=
\cap_{l=1}^{\infty}I_{(w_{k})^{l}}.$  It is easy to check that
the mapping $Q$ is a surjective mapping and there are at most two points in
$P(\sigma_{f}^{*})$ being mapped to a same point under $Q$.

\vskip5pt
{\sc Proposition 3.} {\em
Suppose $f:M\mapsto M$ is a Markov mapping.  Furthermore
suppose
there is the signed scaling function $s_{f}: \Sigma_{f}^{*}\mapsto {\bf R}^{1}$
of $f$.
Then for an $a^{*}$ in $P(\sigma_{f}^{*})$ and
$p=Q(a^{*})$ in $P(f)$, the inverse of the eigenvalue
$e_{f}(p)=(f^{\circ k})'(p)$ of $f$ at $p$ can be calculated by
\[  \frac{1}{e_{f}(p)}=
\prod_{l=0}^{k-1} s_{f}((\sigma_{f}^{*})^{\circ l}(a^{*})).\]}

{\it Proof.} The scale $s((w_{k})^{l})$ at $(w_{k})^{l}$ equals

\[ \frac{
sign(I_{(w_{k})^{l}})I_{(w_{k})^{l}}}
{sign(I_{(w_{k})^{(l-1)} r_{i_{k-1}}\cdots
r_{i_{1}}})I_{(w_{k})^{(l-1)} r_{i_{k-1}}\cdots r_{i_{1}}}}.\]
By using
$ f^{\circ k}(I_{(w_{k})^{l}})=
I_{(w_{k})^{l-1}}$ and the mean value theorem,

\[
s((w_{k})^{l})=
\frac{1}{(f^{\circ k})'(\xi_{(w_{k})^{l}})}\cdot
\frac{ sign(I_{(w_{k})^{l-1}})I_{(w_{k})^{l-1}}}
{ sign(I_{(w_{k})^{l-1}r_{i_{k-1}}\cdots
r_{i_{1}}})I_{(w_{k})^{l-1} r_{i_{k-1}}\cdots r_{i_{1}}}}\]

\[ =\frac{1}{(f^{\circ k})'(\xi_{(w_{k})^{l}})}\cdot
\frac{1}{s((w_{k})^{l-1}r_{i_{k-1}}) \cdots
s((w_{k})^{l-1}r_{i_{k-1}}\cdots r_{i_{1}})} \]
where $ \xi_{(w_{k})^{l}}\in
I_{(w_{k})^{l}}.$

\vskip5pt
Because the maximum $\lambda_{kl}$ of the lengths of the intervals in
$\eta_{kl}$ tends to $0 $ as $l$ goes to infinity, we
have
that $\xi_{(w_{k})^{l}}$ tends to $p$ as
$l$ goes to infinity. Let $l$ tends to infinity, we get

\[ s(a^{*})
= \frac{1}{(f^{\circ k})'(p)} \cdot
\frac{1}{s((w_{k})^{\infty }r_{i_{k-1}})
\cdots s((w_{k})^{\infty} r_{i_{k-1}}\cdots
r_{i_{1}})}\]

\[ =
\frac{1}{(f^{\circ k})'(p)} \cdot
\frac{1}{\prod_{l=1}^{k-1}
s((\sigma_{f}^{*})^{\circ l}(a^{*}))}.\]
This implies Proposition 3.

\vskip5pt
Suppose $f:M\mapsto M$ is a Markov mapping and $\{ \eta_{n}
\}_{n=1}^{\infty}$ is the induced sequence of nested partitions of $M$
by $f$. We say the restriction of $f$ to an interval in $\eta_{1}$ is a
$C^{1+\alpha}$-embedding for some $0< \alpha \leq 1$ if this restriction
and its inverse are
$C^{1}$ with $\alpha$-H\"older continuous derivatives. We say the
$n^{th}$-partition of $M$ by $f$ $
\eta_{n}$ goes to zero exponentially with $n$ if there are positive
constants
$K$ and $\mu <1$ such that $\lambda_{n} \leq K\mu^{n}$ for every
integer $n>0$.

\vskip5pt
{\sc Definition 4.} {\em Suppose $f: M\mapsto M$ is a Markov mapping. We
say $f$ is a good Markov mapping if

\vskip3pt
$(a)$ the restriction of $f$ to every interval in the first partition
$\eta_{1}$ is a $C^{1+\alpha }$-embedding
for some $0< \alpha \leq 1$,

\vskip3pt
$(b)$
the $n^{th}$-partition of $M$ by $f$
$\eta_{n}$ goes to zero exponentially with $n$.}

\vskip5pt
Suppose $f:M\mapsto M$ is a Markov mapping and $\Sigma^{*}_{f}$ is the
dual Cantor set of $f$.  We say a function $s:
\Sigma^{*}_{f}\mapsto {\bf R}^{1}$ is H\"older continuous if there
are positive constants $K$ and $0< \mu <1$ such that
\[ |s(a^{*}_{1})-s(a^{*}_{2})|\leq K\mu^{n} \]
whenever the first $n$ digits of $a^{*}_{1}$ and $a^{*}_{2}$ in $\Sigma_{f}^{*}$
are the same.

\vskip5pt
{\sc Proposition 4.} {\em
Suppose $f:M\mapsto M$ is a good Markov mapping. Then the
signed scaling function
$s_{f}:
\Sigma^{*}_{f}\mapsto {\bf R}^{1}$ exits and is H\"older continuous.}

\vskip5pt
{\it Proof.} The proof of this proposition is the use of the naive
distortion lemma (see [J1] or [J2]) and is similar to the proof of Theorem A in
\S 3. We outline the proof as follows.

\vskip5pt
Suppose $a^{*}_{1}=\cdots w_{n}$ and $a_{2}^{*}=\cdots w_{n}$
are two points in $\Sigma_{f}^{*}$ with the same first $n$ symbols $w_{n}$
(from the right).
Then following the proof of Theorem A and using the naive distortion
lemma (see [J1] or [J2]),
 \[ | s_{f}(a^{*}_{1})-s_{f}(a_{2}^{*})| \leq K |I_{w_{n}}|^{\alpha}\]
where $K$ is a positive constant.
It implies that $s_{f}$ is H\"older continuous.

\vskip5pt
{\sc Example 1.} Suppose
$M$ is the unit interval $[0,1]$ and $l_{0}$, $l_{1}$ and $l_{2}$ are
positive numbers satisfying that $l_{0} + l_{1} +l_{2}=1$.  Let
$K_{0} =[0, l_{0}]$, $K_{1}=[l_{0}, l_{0}+l_{1}]$ and $K_{2}=[l_{0}+l_{1},
1]$
are the subintervals in $M$ and define a Markov mapping $f: M\mapsto M$ as
\[f(x) = \left\{ \begin{array}{ll}
		   \frac{l_{1}+l_{2}}{l_{0}} x + l_{0} &  x \in K_{0},\\
		   -\frac{1}{l_{1}}(x-l_{0}) + 1	 &  x \in K_{1},\\
		   \frac{l_{0} +l_{1}}{l_{2}}(x-l_{0} -l_{1}) & x \in K_{2}.
		  \end{array}
	 \right. \]
\begin{center}
\begin{picture}(120,180)(0,0)
\put (0,40){\line(0,1){120}}
\put (40,40){\line(0,1){120}}
\put (80,40){\line(0,1){120}}
\put (120,40){\line(0,1){120}}
\put (0,40){\line(1,0){120}}
\put (0,80){\line(1,0){120}}
\put (0,120){\line(1,0){120}}
\put (0,160){\line(1,0){120}}
\put (0,80){\line(1,2){40}}
\put (80,40){\line(-1,3){40}}
\put (80,40){\line(1,2){40}}
\put (17,30){$K_{0}$}
\put (57,30){$K_{1}$}
\put (97,30){$K_{2}$}

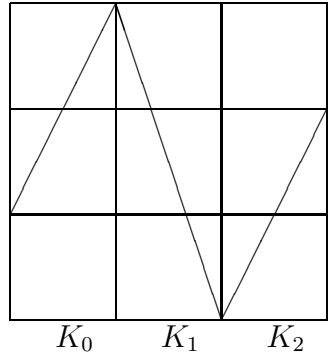
\put (30,15){The graph of $f$}

\put (40,0){Figure 1}
\end{picture}
\end{center}

Suppose $A$ is the induced matrix by the Markov mapping $f$ (see [B]). Then
\[ A = (a_{ij})_{3 \times 3} = \left( \begin{array}{ccc}
						 0 & 1 & 1 \\
						 1 & 1 & 1  \\
						 1 & 1 & 0
						\end{array}
				\right).\]
Then $\Sigma_{f}^{*}$ is
$\Sigma_{A}= \{ a^{*}= (\cdots r_{i_{1}}r_{i_{0}}) |
a_{i_{k}i_{k-1}}=1$ for all $k=1$, $\cdots $, $\infty \}$ and the
scaling function $s_{f}$ is
\[ s_{f}(w) = \left\{ \begin{array}{cl}
			 l_{0} &  a^{*}=( \cdots r_{i_{2}}-1+0),\\
			 \frac{l_{0}}{l_{0}+l_{1}} & a^{*}=( \cdots r_{i_{2}}+2+0),\\
			 -\frac{l_{1}}{l_{1}+l_{2}} & a^{*}=(\cdots r_{i_{2}}+0-1 ),\\
			 -l_{1} & a^{*}=( \cdots r_{i_{2}}-1-1),\\
			 -\frac{l_{1}}{l_{0}+l_{1}} & a^{*}=(\cdots r_{i_{2}}+2-1),\\
			 \frac{l_{2}}{l_{1}+l_{2}} & a^{*}=( \cdots r_{i_{2}}+0+2),\\
			 l_{2} & a^{*} =(\cdots r_{i_{2}} -1+2).
			\end{array}
		\right. \]

\vskip10pt
\centerline {\Large \bf \S 3 The Scaling Function of
a Geometrically Finite}

\centerline {\Large \bf One-dimensional Mapping}

\vskip5pt
A geometrically finite one-dimensional mapping is defined in the
paper [J3]. Let me review this definition here.

\vskip5pt
Suppose $M$ is an oriented connected compact one-dimensional
$C^{2}$-Riemannian
manifold with Riemannian metric $dx^{2}$ and associated length element
$dx$.  Suppose $f: M\mapsto M$
is a $C^{1}$-mapping. We say a point $c\in M$ is a critical point if
the derivative of $f$ at this point is zero.
We say a
critical point $c$ of $f$ is a power law critical
point if
there is a $\gamma \geq 1$ such that
\[ \lim_{x\mapsto c+} \frac{f'(x)}{|x-c|^{\gamma -1}} \hskip5pt and \hskip5pt
\lim_{x\mapsto c-} \frac{f'(x)}{|x-c|^{\gamma-1}} \]
have nonzero limits $A$ and $B$.
We call the numbers $\gamma$ and $\tau =A/B$ the exponent and the
asymmetry of $f$
at the power law critical point $c$ (see [J1] and [J2]). We say a
critical point $c$ of $f$ is
critically finite if the orbit $\{ c$, $f(c)$, $\cdots \}$ is a finite
set.

\vskip5pt
Remember that
an object associated with $f$ is a $C^{1}$-invariant of $f$ if
it is the same for $f$ and for $h\circ f\circ h^{-1}$ whenever $h$ is an
orientation-preserving $C^{1}$-diffeomorphism. We have the
following proposition.

\vskip5pt
{\sc Proposition 5.} {\em Suppose $f:M\mapsto M$ is a $C^{1}$-mapping
and $c$ is a power law critical point of $f$. Then the exponent $\gamma$
and the asymmetry $\tau$ of $f$ at $c$ are $C^{1}$-invariants of $f$.}

\vskip5pt
Henceforth, without loss generality, we will assume that $f$ maps the
boundary of $M$ (if it is not empty) into itself and the one-sided
derivatives of $f$ at all boundary points of $M$ are not zero. We note
that in the general case, a boundary point of $M$ should count as a
critical point anyway.

\vskip5pt
Suppose $f$ has only power law critical points.
We use $CP=\{ c_{1}$, $\cdots $,
$c_{d}\}$ to denote the set of
critical points of $f$ and use
$\Gamma=\{ \gamma_{1}$, $\cdots $, $\gamma_{d}\}$
to denote the
corresponding exponents of $f$. Suppose $\eta_{0}$ is the set of the
intervals in the
complement of the set $CP$ of critical points of $f$ in $M$.

\vskip5pt
{\sc Definition 5.}
{ \em We say the mapping $f$ is $C^{1+\alpha}$ if

\vskip3pt
$(*)$ the restrictions of $f$ to the intervals in
$\eta_{0}$ are $C^{1}$ with $\alpha $-H\"older continuous derivatives
and

\vskip3pt
$(**)$ for every critical point $c_{i}$ of $f$ , there
is a small neighborhood $U_{i}$ of $c_{i}$ in $M$ such that
$\delta_{-,i}(x)=f'(x)/|x-c|^{\gamma_{i}-1}$ for $x<c$ in $U_{i}$ and
$\delta_{+,i}(x)=f'(x)/|x-c|^{\gamma_{i}-1}$ for $x>c$ in $U_{i}$ are
$\alpha$-H\"older continuous.}

\vskip5pt
Suppose the set of the critical orbits $\cup_{n=0}^{\infty}f^{\circ
n}(CP)$ is finite. Then the set of the closures of the intervals
of the complement of the critical orbits $\cup_{n=0}^{\infty}f^{\circ
n}(CP)$ is a Markov partition of $M$ by $f$. We always take this Markov
partition of $M$ by $f$ as the first partition $\eta_{1}$ of $M$ by $f$
in this case.  Suppose $\lambda_{n}$ is the maximum of lengths
of the intervals in the $n^{th}$-partition $\eta_{n}$ of $M$ by $f$.
Remember that the $n^{th}$-partition $\eta_{n}$ tends to zero
exponentially with $n$ if
there are constants $K>0$ and $0< \mu <1$
such that $\lambda_{n}\leq K\mu^{n}$ for
any $n$. The definition of a geometrically finite one-dimensional
mapping is the following.

\vskip5pt
{\sc Definition 6.} {\em  A $C^{1}$-mapping $f:M\mapsto M$ with only
power law
critical points is $C^{1+\alpha}$-geometrically finite for some $0<
\alpha \leq 1$ (geometrically finite) if it satisfies the following
conditions:

\vskip3pt
Smooth condition: $f$ is $C^{1+\alpha }$.

\vskip3pt
Finite condition:
the set of critical orbits $\cup_{i=0}^{\infty }f^{\circ }(CP)\neq
\emptyset $ is finite.

\vskip3pt
No cycle condition: no critical point is periodic.

\vskip3pt
Exponential decay condition: the $n^{th}$-partition
$\eta_{n}$ tends to zero exponentially with $n$.}

\vskip5pt
A technical lemma, the $C^{1+\alpha}$-Denjoy-Koebe
distortion lemma, has been developed in [J2] to study the certain mappings
with finitely many
nonrecurrent critical points. For a geometrically finite one-dimensional
mapping $f: M \mapsto M$, this lemma can be written in the following
simple form (see
of \S 3.3 in [J2]).

\vskip5pt
{\sc Lemma 2} (The $C^{1+\alpha}$-Denjoy-Koebe distortion lemma)
{\em Suppose $f:M\mapsto M$ is geometrically finite.
There are two
positive
constants $A$ and $B$ and a positive integer $n_{0}$ such that for
any inverse branch $g_{n}$ of $f^{\circ n}$ and
any pair $x$ and $y$ in the intersection of one of the intervals in
$\eta_{n_{0}}$ and the domain of $g_{n}$,
the distortion $|g_{n}(x)/g_{n}(y)|$ of $g_{n}$ at these two points
satisfies
\[ \frac{|g_{n}(x)|}{|g_{n}(y)|} \leq \exp
\Big( \Big( A+\frac{B}{D_{xy}} \Big) |x-y|^{\alpha} \Big) \]
where $D_{xy}$ is the distance between $\{ x$, $y\}$ and the
boundary of the domain of $g_{n}$.}

\vskip5pt
\noindent {\bf \S 3.1 The existence of the signed scaling function.}

\vskip5pt
One of the main results in this paper , which is an application of
Lemma 2 (the $C^{1+\alpha }$-Denjoy-Koebe distortion lemma), is the
following:

\vskip5pt
{\sc Theorem A.} {\em
Suppose $f$ is a geometrically finite
one-dimensional mapping. Then
there is the signed scaling function $s_{f}:
\Sigma^{*}_{f}\mapsto {\bf R}^{1}$ of $f$.}

\vskip5pt
{ \it Proof.}
Suppose $U_{1}$, $\cdots $, $U_{d}$ are the
neighborhoods
of the critical points $c_{1}$, $\cdots $, $c_{d}$ of $f$ in
Definition 5. We say an interval $I$ in
$\eta_{n}$ is a critical interval if one of its endpoints is a
critical point of $f$. Suppose
$n_{0}$ is the integer in Lemma 2. Let $n_{1}>n_{0}$ be an integer such
that
every critical interval $I$ in $\eta_{n_{1}}$ is contained in one of
$U_{1}$, $\cdots$, $U_{d}$ and one of its endpoints is not in the
critical orbits
$\cup_{n=0}^{\infty} f^{\circ n}(CP)$. Let ${\cal U}$ be the
union of the critical intervals in
$\eta_{n_{1}}$ and ${\cal V}$ be the union of the non-critical
intervals in $\eta_{n_{1}}$.

\vskip5pt
For a point $a^{*} =\cdots w_{n}$ in $\Sigma_{f}^{*}$, let
$\sigma_{f}^{*}(a^{*}) =\cdots v_{n-1}$ where $w_{n}=v_{n-1}r_{i_{0}}$.
Suppose $I_{w_{n}}$ and $I_{v_{n-1}}$ are the intervals with names $w_{n}$
and $v_{n-1}$, respectively. We note that $I_{w_{n}}$ is a subinterval
of $I_{v_{n-1}}$.
We discuss the sequence $\{
I_{v_{n-1}}\}_{n=1}^{\infty}$ in the two cases. One is that there is
a positive integer $N$ such that $I_{v_{n-1}}$ is contained in ${\cal V}$
for every $n>N$. The other is that there is an increasing subsequence $\{ n_{k}
\}_{k=2}^{\infty}$ of the integers such that $I_{v_{n_{k}-1}}$ is contained
in ${\cal U}$ for every $k\geq 2$. Suppose $n_{2}\geq n_{1}$.

\vskip5pt
In the first case, we use the naive distortion lemma (see [J1] and [J2]) to obtain
the following estimate:

\vskip5pt
For any integers $m>n>N$, the intervals $I_{w_{n}}$ and $I_{v_{n-1}}$
are the images of $I_{w_{m}}$ and $I_{v_{m-1}}$ under $f^{\circ (m-n)}$
and the intervals $I_{w_{N}}$ and $I_{v_{N-1}}$ are the images of
$I_{w_{n}}$ and $I_{v_{n-1}}$ under $f^{\circ (n-N)}$. We note that the
signs of $s(w_{m})$ and $s(w_{n})$ are the same. These imply that
\[ |s(w_{m})-s(w_{n})| = |\frac{ f^{\circ (m-n)}(\xi_{1})}{ f^{\circ
(m-n)}(\xi_{2})} -1| \cdot \frac{ |I_{w_{n}}|}{ |I_{v_{n-1}}|}\]

\[ = |\frac{ (f^{\circ (m-n)})'(\xi_{1})}{ (f^{\circ
(m-n)})'(\xi_{2})} -1| \cdot
|\frac{ (f^{\circ (n-N)})'(\xi_{3})}{ (f^{\circ
(n-N)})'(\xi_{4})}| \cdot \frac{ |I_{w_{N}}|}{ |I_{v_{N-1}}|} \]
which is less than $K |I_{w_{n}}|^{\alpha }$ for a positive constant $K$
where $\xi_{1}$ and $\xi_{2}$ are two points in
$I_{v_{m-1}}$ and $\xi_{3}$ and $\xi_{4}$ are two points in $I_{v_{n-1}}$.
This estimate says that the sequence $\{ s(w_{n}) \}_{n=1}^{\infty}$ is
a Cauchy sequence and the limit
$s(a^{*})=\lim_{n\mapsto \infty} s(w_{n})$ exists.

\vskip5pt
In the other case, we have that the intervals
$I_{w_{n_{k}}}$ and $I_{v_{n_{k}-1}}$
are the images of $I_{w_{n}}$ and $I_{v_{n-1}}$ under $f^{\circ (n-n_{k})}$
for any $n>n_{k}$ and the intervals $I_{w_{n_{2}}}$ and $I_{v_{n_{2}-1}}$ are the images of
$I_{w_{n_{k}}}$ and $I_{v_{n_{k}-1}}$ under $f^{\circ (n_{k}-n_{2})}$. Then
we get

\vskip5pt
\[ |s(w_{n})-s(w_{n_{k}})| = |\frac{ f^{\circ (n-n_{k})}(\xi_{1})}{ f^{\circ
(n-n_{k})}(\xi_{2})} -1| \cdot \frac{ |I_{w_{n_{k}}}|}{
|I_{v_{n_{k}-1}}|}\]

\[ = |\frac{ (f^{\circ (n-n_{k})})'(\xi_{1})}{ (f^{\circ
(n-n_{k})})'(\xi_{2})} -1| \cdot
|\frac{ (f^{\circ (n_{k}-n_{2})})'(\xi_{3})}{(f^{\circ
(n_{k}-n_{2})})'(\xi_{4})}| \cdot \frac{ |I_{w_{n_{2}}}|}{
|I_{v_{n_{2}-1}}|}\]
where
$\xi_{1}$ and $\xi_{2}$ are two points in
$I_{v_{n-1}}$ and
$\xi_{3}$ and $\xi_{4}$ are two points in
$I_{v_{n_{k}-1}}$.
Suppose $L$ is the minimum of lengths of the intervals in
$\eta_{n_{1}}$. Then $D_{\xi_{1}\xi_{2}}$ and $D_{\xi_{3}\xi_{4}}$ in
Lemma 2 (the $C^{1+\alpha}$-Denjoy-Koebe distortion lemma),
respectively,	are both
greater than or equal to $L$ (we can always reduce to this situation).
By using
Lemma 2, there is a positive constant
$K$ such that $|s(w_{n})-s(w_{n_{k}})|\leq K |I_{w_{n_{k}}}|^{\alpha }.$
From this estimate,
$ |s(w_{n})-s(w_{m})|\leq 2K |I_{w_{n_{k}}}|^{\alpha }$ for any $m>n\geq
n_{k}$.
Again the sequence $\{ s(w_{n}) \}_{n=1}^{\infty}$ is
a Cauchy sequence and the limit
$s(a^{*})=\lim_{n\mapsto \infty} s(w_{n})$ exists. We proved Theorem A.

\vskip5pt
{\sc Corollary A1.} {\em Suppose $f$ is a geometrically finite
one-dimensional mapping.
Then the scaling function $|s_{f}|:
\Sigma_{f}^{*}\mapsto {\bf R}^{1}$ of $f$ is
discontinuous.}

\vskip5pt
{\it Proof.}
Suppose $p_{i}$ is the periodic point such
that the critical point $c_{i}$ of $f$ lands on it under some
iterates of $f$ for every $i=1$, $\cdots$, $d_{1}$. Suppose
$OB(p_{i})=\cup_{k=0}^{\infty}f^{\circ k}(p_{i})$ be the periodic orbit of
$p_{i}$ under $f$.  Let $O=\cup_{k=1}^{d_{1}} OB(p_{i})$ be the union of
the periodic orbits $OB(p_{i})$. It is contained in $P(f)$. Suppose $Q:
P(\sigma_{f}^{*})\mapsto P(f)$ is the mapping in Proposition 2. Let ${\cal
A}_{0}$ be the preimage of $O$ under $Q$ and ${\cal A}$ be the union of the
preimages of ${\cal A}_{0}$ under the $n^{th}$-iterate of $\sigma_{f}^{*}$
for $n=0$, $1$, $\cdots $. We claim that all the points in ${\cal A}$ are
discontinuous points of the scaling function $s_{f}: \Sigma_{f}^{*} \mapsto
{\bf R}^{1}$.

\vskip5pt
For an $a^{*}_{0}=\cdots w_{n} \in {\cal A}$, let $I_{w_{n}}$ be the interval
in $\eta_{n}$ with the name $w_{n}$. There is a subsequence $\{ n_{k}
\}_{k=2}^{\infty}$ of the integers such that $I_{w_{n_{k}}}$ tends to a
periodic point
$p_{i}$ of $f$ as $k$ goes to infinity. To simplify our arguments, let
us assume that the critical point $c_{i}$ is not in the
post-critical orbits $\cup_{k=1}^{\infty}f^{\circ
k}(CP)$ and $f(c_{i})=p_{i}$. For a more general situation, the proof
can be easily obtained by modifying the following arguments.

\vskip5pt
Suppose $I_{u_{n_{k}+1}}$ is an inverse branch of the interval
$I_{w_{n_{k}}}$ under $f$ and is contained in
a critical interval $I$ in $\eta_{n_{1}}$ (where $n_{1}$ is the integer in
the proof of Theorem A).  Because the restriction of $f$ to $I$ is
comparable with $|x-c_{i}|^{\gamma_{i}} +f(c_{i})$ for some $0\leq i\leq
l$, then we can show that
the scale $s(u_{n_{k}+1})$ at $u_{n_{k}+1}$ can be calculated as
follows:
\[ s(u_{n_{k}+1})= \frac{ \xi_{n_{k},1}}{\xi_{n_{k},2}} s(w_{n_{k}}). \]
The limit $\mu$ of the sequence $\{ \mu_{k} =\xi_{n_{k}, 1}/\xi_{n_{k}, 2}
\}_{k=2}^{\infty}$ exists and does not equal one if $I_{w_{n_{k}}}\neq
I_{v_{n_{k}-1}}$ where $w_{n_{k}}=v_{n_{k}-1}r_{i_{0}}$ (this is true in
general).  Let $k$ go to infinity, we get that

\[ \lim_{a^{*}= \cdots u_{n_{k}+1} , k\mapsto \infty} s_{f}(a^{*}) =\mu
\cdot s_{f}(a_{0}^{*}).\]
This implies that $s_{f}$ is discontinuous at the point $a_{0}^{*}$.

\vskip5pt
{\bf Remark.} We actually can find all the continuous points and
discontinuous points of $s_{f}$. Let us do it in a little simple case.
Suppose the set $CP$ of critical points of $f$ is
disjoint with the post-critical orbits
$\cup_{n=1}^{\infty} f^{\circ n}(CP)$.
Let ${\cal U}$ and ${\cal V}$ be the sets in the proof of
Theorem A. Suppose $a^{*}=\cdots v_{n-1}r_{i_{0}}$ is a point in the
dual
space $\Sigma_{f}^{*}$. Let $I_{v_{n-1}}$ be the interval in
$\eta_{n-1}$ with
the name $v_{n-1}$.  We say $a^{*}$ is recurrent if there is a subsequence
$\{ n_{k}\}_{k=2}^{\infty}$
of the integers such that $I_{v_{n_{k}-1}}$ is contained in ${\cal U}$ for every
$k\geq 2$. We say $a^{*}$ is totally nonrecurrent if there is an integer $N>0$
such that the preimage of $I_{v_{N-1}}$ under $f^{\circ k}$ is
contained in ${\cal V}$
for every $k\geq 0$. We say $a^{*}$ is wandering if $(a)$ there is an integer $N>0$
such that $I_{v_{n-1}}$ is contained in ${\cal V}$ for every $n\geq N$ and $(b)$
for every $k> N$ there is an integer $n_{k}$
satisfying that
the preimage of $I_{v_{k-1}}$ under $f^{\circ n_{k}}$ intersects with the
interior of ${\cal U}$. Then we can prove that $s_{f}$ is
continuous at all the recurrent and totally nonrecurrent points and
discontinuous at all the wandering points.

\vskip5pt
{\sc Corollary A2.} {\em Suppose $f$ is a geometrically finite
one-dimensional mapping. Then the exponent $\gamma$ of $f$ at a
power law critical point $c$ can be calculated by
the scaling function $s_{f}: \Sigma_{f}^{*}\mapsto {\bf R}^{1}$.}

\vskip5pt
{\it Proof.} Suppose $c_{i_{1}}$, $\cdots $, $c_{i_{n}}$ are critical points
of $f$. We say they form a chain if there are integers
$l_{1}$, $\cdots$, $l_{n-1}$ such that $f^{\circ l}(c_{i_{k}})$ is not a
critical point of $f$ for $ 0< l< l_{k}$ and
$f^{\circ l_{k}}(c_{i_{k}})=c_{i_{k+1}}$.

\vskip5pt
Suppose $c_{i_{1}}$, $\cdots $, $c_{i_{n}}$ form a maximum chain.
Let $I_{w_{m}}$ is an interval in $\eta_{m}$ which has $c_{i_{1}}$ as an
endpoint.  Then $I_{w_{m-m_{k}}}= f^{\circ m_{k}}(I_{w_{m}})$ has
$c_{i_{k}}$ as an endpoint where $m_{k}=l_{1}+\cdots +l_{k}$ for $1\leq
k<n$. Suppose
$l_{n}$ is the smallest integer such that
$p= f^{\circ l_{n}}(c_{i_{n}}) $ is a periodic point of
$f$ and
$I_{w_{m-m_{n}}}=f^{\circ m_{n}}(I_{w_{m}})$ is an interval which has $p$ as
an endpoint where
$m_{n}=l_{1}+\cdots +l_{n}$. Suppose $\gamma_{i_{1}}$, $\cdots $,
$\gamma_{i_{n}}$ are the corresponding exponents of these critical points
and $a^{*}=(w_{m-m_{n}})^{\infty}\in \Sigma_{f}^{*}$. Then
we have that
\[ \gamma_{n} =\frac{ \log |s_{f}(a^{*})|}{\lim_{m\mapsto \infty} \log
|s(w_{m-m_{n-1}})|},\]
and
\[ \gamma_{k} =\frac{ \lim_{m\mapsto \infty} \log
|s(w_{m-m_{k}})|}{\lim_{m\mapsto \infty} \log |s(w_{m-m_{k-1}})|}\]
for $1\leq k< n$.

\vskip10pt
\noindent {\bf \S 3.2 The $C^{1+\alpha }$-classification.}

\vskip3pt
Suppose $f$ and $g$ are geometrically finite and topologically conjugate
by an orientation-preserving homeomorphism $h$. We say $f$ and $g$ are
$C^{1}$-conjugate if $h$ is a
$C^{1}$-diffeomorphism.
One of the corollaries of Proposition 1 and Proposition 5 is that the
signed scaling functions of $f$ and $g$ and the asymmetries of $f$ and
$g$ at
the corresponding critical points are the same if
$f$ and $g$ are $C^{1}$-conjugate.
Another main result in this paper is that the signed
scaling function and the asymmetries at
critical points of a geometrically finite one-dimensional mapping
form a complete set
of $C^{1 }$-invariants within a
topological conjugacy class as follows.

\vskip5pt
{\sc Theorem B.} {\em Suppose $f$ and $g$ are geometrically
finite and
topologically conjugate by an orientation-preserving homeomorphism $h$.
Then $f$ and $g$
are $C^{1}$-conjugate if and only if
they have the same signed scaling function and the same asymmetries at
the corresponding critical points.}

\vskip5pt
{\bf Remark.}
The topological conjugacy
class $[f]$ is the subset of geometrically finite one-dimensional
mappings which are topologically conjugate to $f$. The class $[f]$
equals the union of $[f]_{+}$ and $[f]_{-}$.
Here $[f]_{+}$ is
the subset of geometrically finite one-dimensional mappings which are
topologically conjugate to $f$ by orientation-preserving
homeomorphisms and
$[f]_{-}$ is
the subset of geometrically finite one-dimensional mappings which are
topologically conjugate to $f$ by orientation-reversing
homeomorphisms. There is a one-to-one corresponding
between $[f]_{+}$ and $[f]_{-}$.

\vskip5pt
Actually, we can show
more
on the smoothness of the conjugating mapping $h$
as follows.

\vskip5pt
{\sc Corollary B1.}
{\em
Suppose $f$ and $g$ are $C^{1+\alpha}$-geometrically finite
one-dimensional mappings for some $0< \alpha \leq1$.
Furthermore suppose they are topologically conjugate by an
orientation-preserving homeomorphism $h$.
If $f$ and $g$ have the same signed scaling function and
the same asymmetries at the corresponding critical points, then
$h$ is a $C^{1+\alpha
}$-diffeomorphism. }

\vskip5pt
As we mentioned in the beginning of this subsection, ``if only'' part of
Theorem B is a corollary of Proposition 1 and 5.
We prove ``if'' part of Theorem B and Corollary B1 by
several lemmas.

\vskip5pt
Suppose $f$ and $g$ are geometrically finite and
topologically conjugate by an orientation-preserving homeomorphism $h$,
that is, $h\circ f=g\circ h$. Furthermore suppose $f$ and $g$ are both
$C^{1+\alpha}$ for some $0< \alpha \leq 1$.
We use $\eta_{n, f}$ to denote the $n^{th}$-partition of $M$ by $f$ and
use $\eta_{n, g}$ to denote the $n^{th}$-partition of $M$ by $g$ for
every integer $n\geq 0$. We note that
the dual space $\Sigma_{f}^{*}$ of
$f:M\mapsto M$ and
the dual space $\Sigma_{g}^{*}$ of
$g:M\mapsto M$ are the same.

\vskip5pt
To present a clear idea and to avoid unnecessary notations,
we prove the following lemmas under the assumption that
the set $PC$ of critical points and the post-critical
orbits $\cup_{n=1}^{\infty}f^{\circ n}(CP)$ of $f$ are
disjoint.
We may also assume that there is an interval $I_{k_{0}}$ in the first
partition $\eta_{1,f}$ such that every interval $I$ in the
first
partition $\eta_{1,f}$ covers $I_{k_{0}}$ eventually under
some iterate of
$f$, that means, there is an integer $n$ such that the image of $I$
under
$f^{\circ n}(I)$ contains $I_{k_{0}}$. Otherwise, we can
divide $M$ into finitely many intervals, each of which consists of some
intervals in the first partition $\eta_{1,f}$, such that
the restrictions
of $f$ to these intervals are geometrically finite and satisfy this
assumption.

\vskip5pt
Suppose $A_{f}$, $B_{f}$ and $n_{0,f}$ are the constants in Lemma 2 for
$f$ and $A_{g}$, $B_{g}$ and $n_{0,g}$ are the constants in Lemma 2 for
$g$. Let $A_{0}$, $B_{0}$ and $n_{0}$ are the maximums of $A_{f}$ and
$A_{g}$, $B_{f}$ and $B_{g}$, $n_{0,f}$ and $n_{0,g}$, respectively.

\vskip5pt
We say an interval $I$ in $\eta_{n,f}$
is a critical interval if one
of its endpoints is a critical points of $f$.  Suppose
$n_{1}>n_{0}$
is an integer such that every critical interval $I$ in
$\eta_{n_{1},f}$
has an endpoint which is not in the critical orbits
$\cup_{n=0}^{\infty}f^{\circ n}(CP)$ of $f$.
Suppose $L_{f}$ is the minimum of lengths of the intervals in
$\eta_{n_{1},f}$ and
$L_{g}$ is the minimum of lengths of the intervals in
$\eta_{n_{1},g}$. Let $L$ be
the minimum of $L_{f}$ and $L_{g}$.
We use ${\cal U}$ to denote the union of the critical intervals in
$\eta_{n_{1},f}$ and use ${\cal V}$
to denote the closure of the complement of ${\cal U}$ in $M$.

\vskip5pt
{\sc Lemma B1.} {\em There is a positive constant $K$ such that for
an interval $I$ in $\eta_{n+n_{1}, f}$, if the image $I_{n}=f^{\circ
n}(I)$ of $I$ under $f^{\circ n}$ is a critical interval in
$\eta_{n_{1},f}$, then
\[ \frac{ |f^{\circ m}(z)|}{|f^{\circ m}(w)|} \leq \exp
\Big( K|f^{\circ n}(z)-f^{\circ n}(w)|^{\alpha} \Big) , \]
\[ \frac{ |g^{\circ m}(h(z))|}{|g^{\circ m}(h(w))|} \leq \exp
\Big( K|g^{n}(h(z))-g^{\circ n}(h(w))|^{\alpha } \Big) \]
for any points $z$ and $w$ in $I$.}

\vskip5pt
{\it Proof.} This lemma is actually a corollary of Lemma 2 (the
$C^{1+\alpha}$-Denjoy-Koebe distortion lemma).

\vskip5pt
{\sc Lemma B2.} {\em There is a positive constant $K$ such that for
an interval $I$ in $\eta_{n+n_{1}, f}$, if the image $I_{i}=f^{\circ
i}(I)$ of $I$ under $f^{\circ i}$ is in ${\cal V}$ for any $0\leq i\leq
n$, then
\[ \frac{ |f^{\circ m}(z)|}{|f^{\circ m}(w)|} \leq \exp
\Big( K|f^{\circ n}(z)-f^{\circ n}(w)|^{\alpha} \Big) , \]
\[ \frac{ |g^{\circ m}(h(z))|}{|g^{\circ m}(h(w))|} \leq \exp
\Big( K|g^{n} \Big( h(z) \Big) -g^{\circ n} \Big( h(w) \Big) |^{\alpha }
\Big) \] for any points $z$ and $w$ in $I$.}

\vskip5pt
{\it Proof.} This lemma is actually a corollary of the
naive distortion lemma in [J1] (see also
[J2]).

\vskip5pt
We say a homeomorphism $q: I\mapsto J$ from an interval $I$ to an
interval $J$ is absolutely continuous if it is non-singular
with respect to the Lebesgue measure $m$, that is, $m(X)=0$ if and only
if
$m(h(X))=0$ for any measurable set $X$ of $I$. For example, if $q$ and
$q^{-1}$ are both Lipschitz continuous, then $q$ is absolute continuous.

\vskip5pt
{\sc Lemma B3.} {\em Suppose $h$ is absolutely continuous and
has a differentiable point $p_{0}$ in $I_{k_{0}}$ with nonzero
derivative. Then the restriction of $h$ to every critical interval in
$\eta_{n_{1},f}$ is
$C^{1+\alpha }$.}

\vskip5pt
{\it Proof.}
Suppose $GPI=\cup_{i=0}^{\infty} \cup_{j=0}^{\infty}f^{-j}(f^{\circ
i}(p_{0}))$ is the grand preimage of $p_{0}$ under $f$. Then $GPI$ is a
dense subset of $M$. By the equation
$h\circ f=g\circ h$ and
the definition of $C^{1+\alpha}$ in this paper, $h$ is differentiable at
every point in $GPI$ and there is a constant $K_{0}>0$ such that $h'(x)>
K_{0}$ for all $x\in GPI$.

\vskip5pt
Suppose $I_{w}$ is a critical interval in $\eta_{n_{1},f}$ and
$a^{*}=\cdots
w_{n}w$ is a point in $\Sigma_{f}^{*}$. Let $I_{w_{n}w}$ be the interval
in $\eta_{n_{1}+n,f}$ with name $w_{n}w$. For any pair $x$ and $y$ in
the intersection of
$I_{w}$ and $GPI$, let $x_{n}$ and $y_{n}$ in $I_{w_{n}w}$ be the preimages of
$x$ and $y$ under $f^{\circ n}$. Using the equation $h\circ f=g\circ h$, we
have that
\[ \frac{h'(x)}{h'(y)} =\prod_{n=0}^{\infty}\frac{|f'(y_{n})|}{|f'(x_{n})|}
\prod_{n=0}^{\infty}\frac{|g'(h(x_{n}))|}{|g'(h(y_{n}))|}.\]
By using Lemma B1,
\[ \prod_{n=0}^{\infty}\frac{|f'(y_{n})|}{|f'(x_{n})|} \leq \exp \Big( K
|x-y|^{\alpha } \Big) \]
and
\[ \prod_{n=0}^{\infty}\frac{|g'(h(x_{n}))|}{|g'(h(y_{n}))|}\leq \exp
\Big( K |h(x)-h(y)|^{\alpha } \Big) . \]
This implies that
\[ \frac{h'(x)}{h'(y)} \leq \exp
\Big( K \Big( |x-y|^{\alpha}+|h(x)-h(y)|^{\alpha } \Big) \Big) .\]
From the last inequality, the restriction of $h'$
to the intersection of $I_{w}$ and $GPI$ is uniformly continuous.  Then
it can be extended to a continuous
function on $I_{w}$. Because the restriction of $h$ to $I_{w}$ is absolutely
continuous, this continuous extension is the derivative of the restriction
of $h$ to $I_{w}$. Using the last inequality again, the restriction
of $h$ to $I_{w}$ is $C^{1+\alpha}$.

\vskip5pt
{\sc Corollary B2.} {\em Suppose $h$ is absolutely continuous.
Then the exponents of $f$ and $g$ at the corresponding critical
points are the same.}

\vskip5pt
{\sc Lemma B4.} {\em Suppose $h$ is absolutely continuous.
Then the restriction of $h$ to every interval in $\eta_{n_{1},f}$ is
$C^{1+\alpha }$.}

\vskip5pt
{\it Proof.} We still use the same notations as that in the proof of Lemma
B3. Suppose $I_{w}$ is an interval in
$\eta_{n_{1},f}$ and $a^{*}=\cdots w_{n}w$ is a point in
$\Sigma_{f}^{*}$.
Let $I_{w_{n}w}$ be the interval in $\eta_{n_{1}+n,f}$ with the name
$w_{n}w$.
Suppose $x$ and $y$ are any pair in $I_{w}$ and $x_{n}$ and $y_{n}$ in
$I_{w_{n}w}$ are the
preimages of $x$ and $y$ under $f^{\circ n}$.
Using the equation $h\circ f=g\circ h$, we
have that
\[ \frac{h'(x)}{h'(y)} =\prod_{n=0}^{\infty}\frac{|f'(y_{n})|}{|f'(x_{n})|}
\prod_{n=0}^{\infty}\frac{|g'(h(x_{n}))|}{|g'(h(y_{n}))|}.\]

\vskip5pt
If $I_{w}$ is a critical interval, then it is Lemma B3. Suppose $I_{w}$
is not a critical interval. We consider the sequence $\{ I_{w_{n}w}
\}_{n=0}^{\infty}$ in the two cases. The first is that all $I_{w_{n}w}$ are
contained in ${\cal V}$ and the other is that there is an integer $n$ such
that $I_{w_{n}w}$ is contained in ${\cal U}$.

\vskip5pt
In the first case,
by using Lemma B2,
\[ \prod_{n=0}^{\infty}\frac{|f'(y_{n})|}{|f'(x_{n})|} \leq \exp \Big( K
|x-y|^{\alpha } \Big) \]
and
\[ \prod_{n=0}^{\infty}\frac{|g'(h(x_{n}))|}{|g'(h(y_{n}))|}\leq \exp
\Big( K |h(x)-h(y)|^{\alpha } \Big) . \]
This implies that
\[ \frac{h'(x)}{h'(y)} \leq \exp
\Big (K \Big( |x-y|^{\alpha}+|h(x)-h(y)|^{\alpha } \Big) \Big) .\]
Then by the same arguments in the proof of Lemma B1, we get that the
restriction of $h$ to $I_{w}$ is $C^{1+\alpha}$.

\vskip5pt
For the other case, let $k$ be the smallest integer such that $I_{w_{k}w}$
is in ${\cal U}$. The product
\[ \prod_{n=0}^{\infty}\frac{|f'(y_{n})|}{|f'(x_{n})|}
\prod_{n=0}^{\infty}\frac{|g'(h(x_{n}))|}{|g'(h(y_{n}))|}\]
can be written in
three products
\[ \prod_{n=0}^{k-1}\frac{|f'(y_{n})|}{|f'(x_{n})|}
\prod_{n=0}^{k-1}\frac{|g'(h(x_{n}))|}{|g'(h(y_{n}))|},\]

\[ \frac{|f'(y_{k})|}{|f'(x_{k})|}
\frac{|g'(h(x_{k}))|}{|g'(h(y_{k}))|},\]

\[ \prod_{n=k+1}^{\infty}\frac{|f'(y_{n})|}{|f'(x_{n})|}
\prod_{n=k+1}^{\infty}\frac{|g'(h(x_{n}))|}{|g'(h(y_{n}))|}.\]
By using Lemma B2
\[ \prod_{n=0}^{k-1}\frac{|f'(y_{n})|}{|f'(x_{n})|} \leq \exp \Big( K
|x-y|^{\alpha } \Big) \]
and
\[ \prod_{n=0}^{k-1}\frac{|g'(h(x_{n}))|}{|g'(h(y_{n}))|}\leq \exp
\Big( K |h(x)-h(y)|^{\alpha } \Big) . \]
By using Lemma B1,
\[ \prod_{n=k+1}^{\infty}\frac{|f'(x_{n})|}{|f'(y_{n})|} \leq \exp
\Big( K |x_{k+1}-y_{k+1}|^{\alpha } \Big) \]
and
\[ \prod_{n=k+1}^{\infty}\frac{|g'(h(y_{n}))|}{|g'(h(x_{n}))|}\leq
\Big( K |h(x_{k+1})-h(y_{k+1})|^{\alpha } \Big) . \]

\vskip5pt
Suppose $I_{w_{k}w}$ is in the critical interval $I$ in $\eta_{n_{1},f}$
which has a critical point $c_{i}$ of $f$ as an endpoint. Let
\[ l_{f}(x)=f'(x)/|x-c_{i}|^{\gamma_{i}-1}\hskip5pt and \hskip5pt
l_{g}(h(x))=g'(h(x))/|h(x)-h(c_{i})|^{\gamma_{i}-1}\]
for $x$ in $I$, where $\gamma_{i}$ is the exponent of $f$ at $c_{i}$.
Then
\[ \frac{|f'(y_{k})|}{|f'(x_{k})|}
\frac{|g'(h(x_{k}))|}{|g'(h(y_{k}))|} =
\frac{l_{f}(y_{k})}{l_{f}(x_{k})}\frac{l_{g}(h(y_{k}))}{l_{g}(h(x_{k}))}
\frac{ (\frac{h(y_{k})-h(c_{i})}{|y_{k}-c_{i}|})^{\gamma_{i}-1}}
{ (\frac{h(x_{k})-h(c_{i})}{|x_{k}-c_{i}|})^{\gamma_{i}-1}}. \]
By the definition of $C^{1+\alpha }$ in this paper and Lemma B3,
the functions $l_{f}$,
$l_{g}$ and $h(x)/|x-c_{i}|$ are $\alpha$-H\"older continuous.
There is a positive constant, we still denote it as $K$, such that
\[ \frac{|f'(x_{k})|}{|f'(y_{k})|}
\frac{|g'(h(y_{k}))|}{|g'(h(x_{k}))|} = \exp \Big( K
\Big( |x_{k}-c_{i}|^{\alpha }+|h(x_{k})-h(c_{i})|^{\alpha } \Big) \Big)
. \]

\vskip5pt
All these estimates and the same arguments as that in the proof of Lemma
B1 say that the restriction of $h$ to $I_{w}$ is $C^{1+\alpha }$.

\vskip5pt
{\sc Lemma B5.} {\em Suppose $h$ is absolutely continuous.
Furthermore suppose $f$ and $g$ have the same asymmetries at the
corresponding critical points. Then
$h$ is $C^{1+\alpha}$.}

\vskip5pt
{\it Proof.}
From Lemma B4, the restriction of $h$ to every interval in
$\eta_{n_{1},f}$ is $C^{1+\alpha
}$. This implies that for every interval $I$ in $\eta_{n_{1},f}$, the
one-sided
limits of the derivative $h'|I$ at the endpoints of $I$
exist.
We need to prove that these one-sided limits are the same at a common
endpoint of any two intervals in $\eta_{n_{1},f}$.

\vskip5pt
Suppose $I$ and $I'$ are two intervals in $\eta_{n_{1},f}$
and have a common endpoint $p$. By the equation $h\circ f=g\circ h$,
\[ h'(p-) = \lim_{x\mapsto p-} \frac{ f'(x)}{ g'(h(x))} h'(p_{1}-), \]
\[ h'(p+) = \lim_{x\mapsto p+} \frac{ f'(x)}{ g'(h(x))} h'(p_{1}+) \]
where $p_{1}$ is a point in the preimage of $p$ under $f$. Without loss
generality, we may assume that $p_{1}$ is an interior point of an interval
in $\eta_{n_{1}}$. Then $h'(p_{1}-)=h'(p_{1}+)$.

\vskip5pt
If $p$ is not a critical point of $f$, it is easy to see that
$h'(p-)=h'(p+)$.

\vskip5pt
Suppose $p$ is a critical point $c_{i}$ of $f$. Let
\[ A_{f}(p)=\lim_{x \mapsto p-} f'(x)/|x-p|^{\gamma_{i}-1},\]

\[ B_{f}(p)=\lim_{x \mapsto p+} f'(x)/|x-p|^{\gamma_{i}-1},\]

and

\[ A_{g}(h(p))=\lim_{x \mapsto p-} g'(h(x))/|h(x)-h(p)|^{\gamma_{i}-1},\]

\[ B_{g}(h(p))=\lim_{x \mapsto p+} g'(h(x))/|h(x)-h(p)|^{\gamma_{i}-1}.\]
Then
\[ (h'(p-))^{\gamma_{i}} = \frac {A_{f}(p)}{A_{g}(h(p))} h'(p_{1}),\hskip5pt
and \hskip5pt
(h'(p+))^{\gamma_{i}} = \frac {B_{f}(p)}{B_{g}(h(p))} h'(p_{1}).\]
The equality
\[ A_{f}(p)/B_{f}(p)=A_{g}(h(p))/B_{g}(h(p))\]
implies that
$h'(p-)=h'(p+)$. We proved Lemma B5.

\vskip5pt
{\sc Lemma B6.} {\em Suppose $\{ a_{n} \}_{n=0}^{\infty}$ and $\{ b_{n}
\}_{n=0}^{\infty}$ are two sequences of positive numbers and there is a
constant $K>0$ such that $a_{n}/b_{n} \leq K$ for any $n$. Then
$(\sum_{n=0}^{\infty} a_{n})/ (\sum_{n=0}^{\infty} b_{n}) \leq
K.$}

\vskip5pt
The proof of this lemma is very easy, but it is very useful in the study of
dynamic systems.

\vskip5pt
{\sc Lemma B7.} {\em Suppose $f$ and $g$ have the same scaling function. Then
the conjugating mapping $h$ is Lipschitz continuous.}

\vskip5pt
{\it Proof.}
Suppose $N$ is the number of the intervals in the first
partition $\eta_{1,f}$ and $n_{1}$ is greater than $2N$.

\vskip5pt
For every
integer $n\geq 0$
and every interval $I_{ww_{n}}$ in $\eta_{n+n_{1}, f}$ with name
$ww_{n}$, let
$s_{f}(ww_{n})$ and $s_{g}(ww_{n})$ be
the scales at $ww_{n}$ (with respect to $f$ and $g$).
Then we have an equation
\[ \frac{
|h(I_{ww_{n}})|}{|I_{ww_{n}}|}=
\frac{|s_{g}(ww_{n})|}{|s_{f}(ww_{n})|}
\frac{|h(I_{wv_{n-1}})|}{|I_{wv_{n-1}}|},\]
where
$ww_{n}=wv_{n-1}r_{i_{0}}$.

\vskip5pt
Let $a^{*}=\cdots u_{m}ww_{n}$ be a point in $\Sigma_{f}^{*}$ and
$I_{u_{m}ww_{n}}$ be the interval in $\eta_{m+n_{1}+n,f}$ with the name
$u_{m}ww_{n}$. We discuss the sequence $\{ I_{u_{m}ww_{n}}
\}_{m=0}^{\infty}$ in the three cases.
The first case is that $I_{ww_{n}}$ is in ${\cal U}$
The second case is that all $I_{u_{m}ww_{n}}$ are
in ${\cal V}$. The third case is that there is a positive integer $k$ such that
$I_{u_{m}ww_{n}}$ is in ${\cal V}$ for every $0 \leq m\leq k$ and
$I_{u_{k+1}ww_{n}}$ is in ${\cal U}$.

\vskip5pt
In the first and the second cases, we use Lemma B2 and Lemma B1,
respectively, to prove the following:

\vskip5pt
There is a constant $0 < \mu<1$ such that
\[ ||s_{f}(a^{*})|-|s_{f}(ww_{n})||\leq \exp (\mu^{n}),\]
\[ ||s_{g}(a^{*})|-|s_{g}(ww_{n})||\leq \exp (\mu^{n}).\]
Because there is a constant $\beta>0$ such that
$|s_{f}(c^{*})|\geq \beta $
for all $c^{*}\in \Sigma_{f}^{*}$ and $|s_{g}|=|s_{f}|$, we can find
a constant, we still denote it as $\mu$, in $(0, 1)$ such that
\[ \frac{ |s_{f}(ww_{n})|}{ |s_{g}(ww_{n})|} \leq \exp (\mu^{n}).\]

\vskip5pt
For the third case, let us suppose that $I_{u_{k+1}ww_{n}}$ is contained in
a critical interval $I$ in $\eta_{n_{1},f}$ which has a
critical point $c$ of $f$
as an endpoint. There is an integer $ 0< m< 2N$ such that $f^{\circ
m}(c)$ is a periodic point of $f$ and the interval
$f^{\circ m}(I_{u_{k+1}ww_{n}})$
is contained in $f^{\circ m}(I)$. We note that
$f^{\circ m}(I)$, which has $p$ as an endpoint, is an interval in
$\eta_{n_{1}-m, f}$ and $n_{1}-m >0$. We also note that
$f^{\circ m}(I_{u_{k+1}ww_{n}})$ is an interval in
$\eta_{n_{1}+n+k+1-m,f}$. Now we can find a point $b^{*}=\cdots v_{j}$
in $\Sigma_{f}^{*}$ such that the first $n_{1}+n+k+1-m$ symbols of
$a^{*}$ and $b^{*}$ (from the right) are the same and the interval
$I_{v_{j}}$ in $\eta_{j}$ is contained in ${\cal V}$ for every $j>0$
and tends to the periodic orbit
$\cup_{i=0}^{\infty}f^{\circ i}(p)$ as $j$ goes to infinity.
By using Lemma B2, there is constant, we still denote it as $\mu$, in
$(0,1)$ such that
\[ ||s_{f}(b^{*})|-|s_{f}(ww_{n})||\leq \exp (\mu^{n}),\]
\[ ||s_{g}(b^{*})|-|s_{g}(ww_{n})||\leq \exp (\mu^{n}).\]
Again, because
$|s_{f}(c^{*})|\geq \beta $
for all $c^{*}\in \Sigma_{f}^{*}$ and $|s_{g}|=|s_{f}|$, we can find a
constant, we still denote it as $\mu$, in $(0, 1)$ such that
\[ \frac{ |s_{f}(ww_{n})|}{ |s_{g}(ww_{n})|} \leq \exp (\mu^{n}).\]

\vskip5pt
Suppose $K_{0}$ is the minimum of the ratios,
$|h(I_{w})|/|I_{w}|$, for $I_{w}$ in $\eta_{n_{1},f}$. By the above
arguments and the induction, we find a
sequence $\{ K_{n} \}_{n=0}^{\infty}$
and a constant $\mu\in (0, 1)$
such that
\[ \frac{
|h(I_{ww_{n}})|}{|I_{ww_{n}}|} \leq K_{n} \]
for every interval $I_{ww_{n}}$ in $\eta_{n+n_{1},f}$, $n \geq 0$,
and \[ K_{n} \leq \exp (\mu^{n}) K_{n-1}\]
for every integer $n\geq 1$.
This yields a positive constant $K$ such that
\[ \frac{
|h(I_{ww_{n}})|}{|I_{ww_{n}}|} \leq K \]
for every $n\geq 0$ and every interval $I_{ww_{n}}$ in
$\eta_{n+n_{1},f}$.

\vskip5pt
Because the union of the boundary points of all the intervals in
$\eta_{n+n_{1},f}$ for all the integer $n \geq
0$ is a dense subset in $M$, by using Lemma B6
\[ \frac{|h(x)-h(y)|}{|x-y|} \leq K \]
for every pair $x$ and $y$ in $M$. In the other words, $h$ is Lipschitz
continuous.

\vskip5pt
{\sc Lemma B8.} {\em Suppose $f$ and $g$ have the same scaling
function and the same asymmetries at the corresponding periodic
points.
Then the conjugating mapping $h$ is a $C^{1+\alpha }$-diffeomorphism.}

\vskip5pt
{\it Proof.} From Lemma B7, the mapping $h$ is Lipschitz
continuous. It is then differentiable at almost every points in
$M$. Let $p_{0}$ be a point in $I_{k_{0}}$ such that $h$ is differentiable
at this
point. Suppose $GPI = \cup_{i=0}^{\infty} \cup_{j=0}^{\infty}
f^{-j}(f^{\circ i}(p_{0}))$ is the grand preimage of $p_{0}$ under $f$. It is a
dense subset of $M$.
If the derivative $h'(p_{0})$ at $p_{0}$ is zero, then by the equation
$h\circ f=g\circ h$, the derivative $h'(p)$ at every point $p\in GPI$ is zero.
But $h$ is absolutely continuous, this implies that $h$ is a constant.
So the derivative $h'(p_{0})$ is not zero. Now Lemma B6 says that $h$
is $C^{1+\alpha}$. The same arguments can be applied to $h^{-1}$.
Hence $h$ is a $C^{1+\alpha }$-diffeomorphism.

\vskip5pt
Lemma B1 to Lemma B8 give the proof of Theorem A.

\vskip5pt
Suppose $f:M\mapsto M$ is a geometrically finite one-dimensional mapping
and
$s_{f}:\Sigma_{f}\mapsto \Sigma_{f}$ is the signed scaling function of
$f$. The eigenvalue $e_{f}(p)=(f^{\circ n})'(p)$
of $f$ at a periodic point $p$ of period $n$ and
the exponent $\gamma$ of $f$ at a critical point $c$
can be calculated by the signed scaling function $s_{f}$ of $f$ showed
by Proposition 2 and
by Corollary A2.  Both of them are then clearly
$C^{1}$-invariants.
In the case that the set of periodic points of $f$ is dense in $M$, we
can show that the eigenvalues of $f$ at
periodic points and the exponents and the asymmetries of $f$ at critical
points form a complete $C^{1}$-invariants within a topologically
conjugate class as
follows.

\vskip5pt
{\sc Theorem C.} {\em Suppose $f$ and $g$ are
$C^{1+\alpha}$-geometrically finite one-dimensional mappings for some
$0< \alpha \leq 1$. Furthermore,
suppose the
set of periodic points of $f$ is dense in $M$ and
suppose $f$ and $g$ are topologically
conjugate by an orientation-preserving homeomorphism $h$
If $f$ and $g$ have the
same eigenvalues at the corresponding
periodic points and the same exponents at the
corresponding critical
points, then they have the same scaling functions.
Moreover, if $f$ and $g$ have also the same asymmetries at the
corresponding critical points, then $h$ is a $C^{1+\alpha
}$-diffeomorphism.}

\vskip5pt
{\it Proof.} The idea of the proof of Theorem C is the same as that of
the proof of Theorem B and that of the proof of Theorem 1.4 in [J1,
p 63-74]. The details will be omitted.

\end{document}